\documentclass[12pt,letterpaper]{article}
\usepackage[latin1]{inputenc}
\usepackage[top=3cm, bottom=3cm, left=3cm, right=3cm]{geometry}
\usepackage{authblk, graphicx,rotating,amsbsy,amssymb,amsmath, hyperref, natbib, bbm, color}

\newtheorem{proposition}{Proposition}
\newtheorem{example}{Example}

\begin{document}

\title{\Large Restricted type II maximum likelihood priors \\ on regression coefficients}

\author[1]{V\'ictor Pe\~na\thanks{Corresponding author, email: victor.pena@baruch.cuny.edu}}
\author[2]{James O. Berger}
\affil[1]{Paul H. Chook Department of Information Systems and Statistics, Baruch College, The City University of New York}
\affil[2]{Department of Statistical Science, Duke University}

\maketitle

\begin{abstract}
In Bayesian hypothesis testing and model selection, prior distributions must be chosen carefully. For example, setting arbitrarily large prior scales for location parameters, which is common practice in estimation problems, can lead to undesirable behavior in testing (Lindley's paradox). 
We study the properties of some restricted type II maximum likelihood (type II ML) priors on regression coefficients. In type II ML, hyperparameters are ``estimated'' by maximizing the marginal
likelihood of a model. In this article, we define priors by estimating their  variances or covariance matrices, adding restrictions which ensure that the resulting priors are at least as vague as conventional proper priors for model uncertainty.
We find that these type II ML priors typically yield results that are close to answers obtained with the Bayesian Information Criterion (BIC; \cite{schwarz1978estimating}).
\end{abstract}

\newpage

\section{Introduction}

In this article, we investigate the properties of restricted type II maximum likelihood (type II ML) priors on regression coefficients under model uncertainty. Along the way, we establish connections with the Bayesian Information Criterion (BIC; \cite{schwarz1978estimating}) and proper priors. Operationally, parametric type II ML proceeds as follows: (1) start with a parametric model for the data $y$, specified with a sampling density $f(y \mid \theta)$ and a prior $\pi_\eta(\theta)$ that depends on a hyperparameter $\eta \in \mathcal{C}$ and (2) set $\eta$ by maximizing the marginal likelihood $m(y)$ of the model, that is
$$
\widehat{\eta} =  \arg  \max_{\eta \in \mathcal{C}} \int f(x \mid \theta) \pi_\eta (\theta) \, \mathrm{d} \theta  = \arg  \max_{\eta \in \mathcal{C}} m(y).
$$
Type II ML was named and extensively studied in \cite{good1965estimation}, and it can be seen as a particular instance of empirical Bayes which, in general, ``estimates'' the hyperparameter $\eta$ from the data (although not necessarily by maximizing the marginal likelihood: a popular alternative is the method of moments).

The motivation for this work was to seek a compromise between the use of conventional priors for model uncertainty (e.g., Zellner's $g$-priors
or Zellner-Siow prior; \cite{zellner1980posterior, zellner1986assessing}) and BIC. Conventional
priors are typically centered at the smallest (or null) model and can be quite far from the likelihood function arising from a larger model. Prior distributions centered at null models can be oriented in directions away from the likelihood function, which would seem to unduly favor the null model.

\cite{raftery1995bayesian} shows that
BIC is a good approximation to the marginal likelihood one obtains when a normal prior that is centered at the maximum likelihood estimate (MLE) $\theta$ is used. This is actually a type II ML prior, arising from estimating the prior mean by type II ML. However, it seems like an extreme use of type II ML because it centers the prior completely around the model likelihood function.

The compromise studied
herein is to keep the prior centered at the null model, as with current conventional priors, but allow the prior variance or covariance matrix
to be estimated by type II ML. We hoped that this would strike a balance between conventional priors and BIC but, while we find that this ``variance-oriented" type II ML prior does yield compromise results, the conclusions are typically closer to BIC.

A second surprise was that the ``variance-oriented" type II ML prior (and resulting Bayes factors) can be computed in closed form, even for the case of entire
unknown prior covariance matrices (this is a computational advantage over, e.g., the Zellner-Siow priors). The
importance of making  restrictions on the hyperparameters is also highlighted; without them, one can even have inconsistent model selection (for example, if the scale parameter $g$ of a Zellner $g$-prior \citep{zellner1986assessing} is estimated without restrictions, the resulting procedure is not consistent if the null model is true \citep{liang2008mixtures}).

Our work is partially motivated by  \cite{bayarri2019prior}, where prior-based versions of BIC (named PBIC and PBIC*) are defined. In particular, PBIC* is a version of BIC which builds upon a restricted type II ML version of the so-called ``robust'' prior \citep{berger1985statistical}. The scales of the prior in PBIC* maximize an approximate marginal likelihood subject to a unit-information restriction. The fact that PBIC* is well-behaved in the examples covered in \cite{bayarri2019prior} motivated us to study the properties of restricted type II ML procedures under model uncertainty in greater detail. 

The scenarios we consider in this article involve regression coefficients in normal linear models (Section~\ref{sec:linmod}),  high-dimensional ANOVA  (Section~\ref{sec:ANOVA}), and the nonparametric regression example in \cite{shibata1983asymptotic} (Section~\ref{sec:shibata}). In this latter
section we also highlight how type II ML can be fruitfully used when prior information is available. The article ends with conclusions. All the proofs are relegated to the supplementary material.

\section{Type II ML priors in normal linear models} \label{sec:linmod}

\subsection{Derivation of the type II ML prior}
Consider the normal linear model
$$
Y  = X_0 \beta_0 + X \beta + \epsilon, \, \, \epsilon \sim N_n(0_n, \sigma^2 I_n),
$$
where $Y \in \mathbb{R}^n$, $X_0 \in \mathbb{R}^{n \times p_0}$ contains common predictors, and $X \in \mathbb{R}^{n \times p}$ contains model-specific predictors. We assume that the predictors are linearly independent and the common and model-specific parameters are orthogonal, so that $X_0'X = 0_{p_0 \times p}$ (if $X_0 = 1_n$, this amounts to centering $X$). In this section, the prior on the common parameters is the right-Haar prior $\pi(\beta_0, \sigma^2) \propto 1/\sigma^2$, which is supported by group invariance arguments in \cite{berger1998invariant} and \cite{bayarri2012criteria}.

The prior distribution we consider for $\beta$, given $\sigma^2$, is the $N_p(\beta \mid 0_p, \sigma^2 W)$ normal prior with mean $0_p$ and
positive definite covariance matrix $W$.
For a fixed $W$ and $n \ge p+p_0$, the marginal likelihood is
\begin{align*}
m_W(Y) &= \int_{\mathbb{R}} \int_{\mathbb{R}^p} \int_{\mathbb{R}_+} N_n(Y \mid X_0 \beta_0 + X \beta, \sigma^2 I_n) \, N_p(\beta \mid 0_p, \sigma^2 W) \, 1/\sigma^2 \, \mathrm{d}  \beta_0 \, \mathrm{d}  \beta \, \mathrm{d}  \sigma^2 \nonumber \\
&=  \frac{\Gamma\left(\frac{n-p_0}{2}\right) \pi^{-(n-p_0)/2}}{\left( |X' X| \, |X_0' X_0| \, |(X'X)^{-1}+W|  \right)^{1/2}} [ \mathsf{SSE} + \widehat{\beta}'[W + (X'X)^{-1}]^{-1} \widehat{\beta} ]^{-\frac{(n-p_0)}{2}}, \label{marginal}
\end{align*}
where $\widehat{\beta} = (X'X)^{-1} X' Y$,
$\mathsf{P}_{X} = X (X' X)^{-1} X'$, $\mathsf{P}_{X_0} = X_0(X_0'X_0)^{-1}X_0'$, and $\mathsf{SSE} = Y'(I_n - \mathsf{P}_{X_0}-\mathsf{P}_{X})Y$.
The type II ML approach to determination of $W$ consists in maximizing the marginal likelihood over $W$, using the result as the prior
covariance matrix. An earlier version of this (see \cite{george2000calibration, hansen2003minimum, liang2008mixtures}) considered
$g$-priors arising from $W$ of the form $W= g \,  \sigma^2 (X'X)^{-1}$, and then maximizing the marginal likelihood over the choice of $g$.

While this maximization over $W$ can be done in closed form, the result is not satisfactory, in that the result is
a singular matrix. We will circumvent this issue by constraining $W$ under the maximization, and will do
so through the concept of a ``unit-information prior."
The expected Fisher information of the regression coefficient $\beta$ is $(X'X)/\sigma^2$, so one can argue that $(X'X)/(n \sigma^2)$ contains as much information as a ``typical'' observation in the sample \citep{kass1995reference, raftery1995bayesian,hoff2009first}. The $N_p(0,n \sigma^2(X'X)^{-1})$ prior is often referred to as the unit-information (normal) prior, and it is a reasonably vague (but necessarily proper) prior for dealing with model uncertainty. Motivated by this discussion, we study the restricted type II ML prior
\begin{align*}
\beta \mid \sigma^2 &\sim N_p(0_p, \sigma^2 \, \widehat{W}) \\
\widehat{W} &= \mathrm{arg \, max}_{W \succeq n(X'X)^{-1}} m_W(Y) \,,
\end{align*}
where $A \succeq B$ means that $A-B$ is positive semidefinite. This will ensure that the restricted type II ML covariance will be
at least as disperse as the unit-information prior covariance. The lower bound is also an instance of Zellner's $g$-prior where $g=n$.

 In the context of estimation,  \citet{dasgupta1988frequentist}, \citet{leamer1978specification}, and \citet{polasek1985sensitivity} study priors that resemble our type II ML prior, bounding the prior covariance matrix both above and below.

Proposition~\ref{prop:wmatrix} below shows that the covariance matrix that maximizes $m_W(Y)$ subject to $W \succeq n(X'X)^{-1}$ is a linear combination of the unrestricted maximum over all positive semidefinite matrices, which is proportional to $\widehat{\beta} \widehat{\beta}'$, and the lower bound $n(X'X)^{-1}$.

\begin{proposition} \label{prop:wmatrix}
For $n > p+p_0$, the solution to the optimization problem
\begin{align*}
&\mathrm{maximize} \, \,  m_W(y) \\  &\mathrm{subject \, to}  \, \,  W \succeq n(X'X)^{-1}
\end{align*}
can be written as
\begin{align*}
	\widehat{W} &= a \,\widehat{\beta} \widehat{\beta}' + n(X'X)^{-1} \\
	a &= \max \{0, (n-p_0-1)/\mathsf{SSE}-(n+1)/\mathsf{SSR}\} \\
	\mathsf{SSR} &= \widehat{\beta}'X'X \widehat{\beta} \,.
\end{align*}
\end{proposition}
In the following subsections, we study the properties of the type II ML prior on $\beta$ that takes $\widehat{W}$ as its covariance matrix in model selection and uncertainty, estimation, and prediction.

\subsection{Model uncertainty and selection} \label{sec:modunc}
Let $X_i$ be a design matrix that includes a subset of $p_i$ out of the $p$ predictors in $X$, with $i \in \{1, 2, \, ... \, , 2^p\}$ ($p_i$ can be 0, which corresponds to the null model), and let $\mathcal{M}_i$ be the model $Y = X_0 \beta_0 + X_i \beta_i + \epsilon_i$, where $\epsilon_i \sim N_n(0,\sigma^2 I_n)$ and $X_0' X_i = 0_{p_0 \times p_i}$. Throughout, we set prior covariance matrices locally -- that is, each $\mathcal{M}_i$ is assigned its own $\widehat{W}_i$. The local approach to empirical Bayes model selection is justified through information-theoretical arguments in \cite{hansen2003minimum}. We perform model selection using null-based Bayes factors, namely
$$
\mathsf{BF}_{i0} =  \frac{\int N_n(Y \mid X_0 \beta_0 + X_i \beta_i, \sigma^2 I_n) \pi_{\mathsf{ML}}(\beta_0, \beta_i, \sigma^2)  \, \mathrm{d}(\beta_0, \beta_i, \sigma^2)}{\int N_n(Y \mid X_0 \beta_0, \sigma^2 I_n) \pi_0(\beta_0, \sigma^2) \, \mathrm{d} (\beta_0, \sigma^2)} = \frac{m_i(Y)}{m_0(Y)}.
$$
We use the notation $\pi_{\mathsf{ML}}$ for the joint (type II ML) prior under $\mathcal{M}_i$ $\pi_{\mathsf{ML}}(\beta_0, \beta_i, \sigma^2) \propto N_{p_i} ( \beta_i \mid 0_{p_i} , \, \sigma^2 \widehat{W}_i) \, 1/\sigma^2 $. The prior under the null model is $\pi_0(\beta_0, \sigma^2) \propto 1/\sigma^2$.
Combining the result in Proposition \ref{prop:wmatrix} with the Sherman-Morrison formula and the matrix determinant lemma (which can be found, for example, as Equations 160 and 24 in \cite{petersen2008matrix}, respectively), it is straightforward to see that the null-based Bayes factor of $\mathcal{M}_i$ under the type II ML covariance matrix is
\begin{equation} \label{eq:bf}
  \mathsf{BF}_{i0} = \begin{cases}  (n+1)^{\frac{n-p_0-p_i}{2}} \left[ n (1-R^2_i) + 1 \right]^{-(n-p_0)/2} &  \text{ if } R_i^2 \le \frac{n+1}{2n-p_0} \\ \varphi(n)^{-1/2} \left( R_i^2\right)^{-1/2} (1-R^2_i)^{-(n-p_0-1)/2} &  \text{ if } R_i^2 > \frac{n+1}{2n-p_0},\end{cases}
\end{equation}
where $\varphi(n) =  [(n+1)^{p_i-1}(n-p_0)^{n-p_0}]/[(n-p_0-1)^{n-p_0-1}]$ and $R_i^2 = 1-\mathsf{SSE}_i/\lVert Y - X_0 \widehat{\beta}_0 \rVert^2$. The first case corresponds to the null-based Bayes factor with the lower bound $\widehat{W}_i = \sigma^2 n(X_i'X_i)^{-1}$. The type II ML procedure differs from the lower bound only if the signal-to-noise ratio (that is, $R^2_i$) is high enough. This feature prevents the procedure from unduly favoring larger models. 

Before we study the properties of the prior in more detail, we present an example with $p = 2$ predictors to introduce some geometric intuition. In addition, the example  will help us highlight that the lower bound prior $\pi_{\mathsf{LB}}(\beta_i \mid \sigma^2) = N_{p_i}(\beta_i \mid 0_{p_i}, \sigma^2 n(X_i'X_i)^{-1})$ has a particular asymmetry with respect to the sign of the correlation between the predictors. It also serves as motivation to compare $\pi_{\mathsf{LB}}$ and the type II ML prior $\pi_{\mathsf{ML}}$ to the Bayesian Information Criterion (BIC; \cite{schwarz1978estimating}), which is defined as
$$
-2 \log N_n(Y \mid X_0 \widehat{\beta}_0 + X_i \widehat{\beta}_i, \widehat{\sigma}^2_i I_n) + p_i \log n,
$$
 where $\widehat{\beta}_0$ $\widehat{\beta}_i$, and $\widehat{\sigma}^2_i$ are the maximum likelihood estimators of $\beta_0, \beta_i$ and $\sigma^2$, respectively. Throughout, we treat $\exp(-\mathsf{BIC}/2)$ as an approximate marginal likelihood with the understanding that the ``BIC'' of the null model is $-2 \log N_n(Y \mid X_0 \widehat{\beta}_0, \widehat{\sigma}^2_0 I_n)$. These choices lead to the null-based Bayes factor
$$
\mathsf{BF}_{i0, \mathsf{BIC}} = n^{-p_i/2} (1-R^2_i)^{-n/2}.
$$
\cite{raftery1995bayesian} observed that $\exp(-\mathsf{BIC}/2)$ is an excellent approximation to the marginal likelihood arising
from $\pi_{\mathsf{BIC}}(\beta_i \mid \sigma^2) = N_{p_i}({\beta}_i \mid \widehat{\beta}_i , \sigma^2 n(X_i'X_i)^{-1})$, which is
Zellner's $g$-prior with $g = n$, but centered at $\widehat{\beta}_i$ instead of $0_{p_i}$. Indeed, under such type II ML prior, the null-based Bayes factor is
$$
\mathsf{BF}_{i0, \widehat{\beta}} = (n+1)^{-p_i/2} (1-R_i^2)^{-(n-p_0)/2},
$$
which is almost identical to $\mathsf{BF}_{i0, \mathsf{BIC}}$. Another prior we will be considering in numerical comparisons is the Zellner-Siow prior, which is $\mathrm{Cauchy}_{p_i}(0_{p_i}, \sigma^2 n(X_i'X_i)^{-1})$, since this is one of the most commonly recommended model uncertainty priors.

\begin{figure}[h!]
  \centering
    \includegraphics[width=\textwidth]{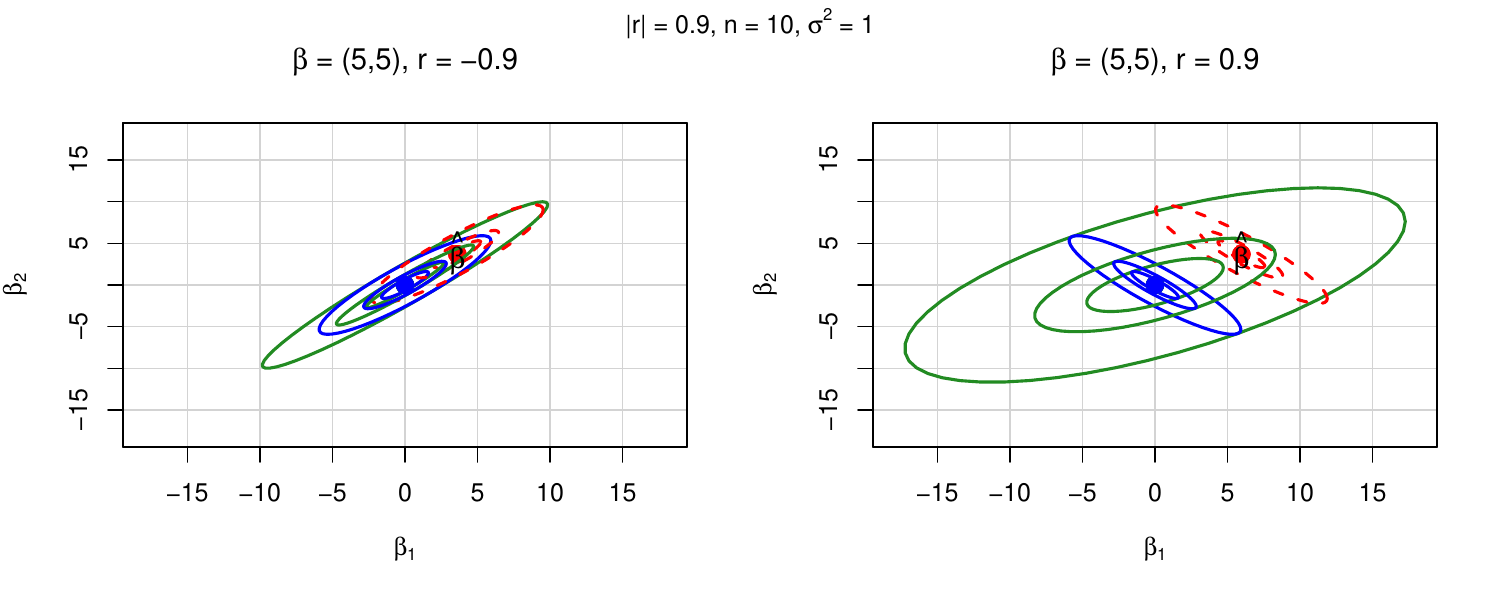}
    \caption{Highest probability density regions (20\%, 50\%, 95\%) of the lower bound (g-prior) $N_p(0_p, n(X'X)^{-1})$ (solid blue), ``BIC prior'' $N_p(\widehat{\beta}, n(X'X)^{-1})$ (dashed red), and the type II ML prior $N_p(0_p, \widehat{W})$ (solid green). The MLE is indicated with a $\widehat{\beta}$ symbol. }
    \label{fig:contours}
\end{figure}

\begin{example} (Correlated predictors) \label{ex:twopreds} Consider a model with 2 standardized (centered and scaled) predictors and an intercept, $Y = 1_n \alpha + X \beta + \epsilon$ where $\beta = (\beta_1, \, \beta_2)'$ and $\epsilon \sim N_n(0_n, I_n)$. Since the predictors are standardized, their (uncorrected) sample correlation is the off-diagonal entry of $(X'X)/n$, which we denote $r$. The prior covariance between $\beta_1$ and $\beta_2$ implied by the prior $\beta \mid \sigma^2 = 1 \sim N_2(0_2, n(X'X)^{-1})$ is $-r/(1-r^2)$ \citep{ghosh2015}. Therefore, if $X_1$ and $X_2$ are positively correlated, the prior covariance between $\beta_1$ and $\beta_2$ induced by the prior is negative (and conversely for negative correlations).

We set $n = 10$, $\beta = (5,5)'$ and consider two cases: $r = 0.9$ and $r = -0.9$. In order to isolate the effect of changing the sign of $r$ as much as possible, we use the same random $\epsilon$ in both cases and the same $N_1(0,1)$ random numbers for generating the design matrices before transforming them (deterministically, via principal component scores times the Choleski matrix square-root of the target sample covariance) to correlated predictors with the desired $r$.

Figure~\ref{fig:contours} shows contours of $N_2(0_2, n(X'X)^{-1})$ (solid blue) and $N_2(0, \widehat{W})$ (solid green; setting $\sigma^2 = 1$), the type II ML prior. It also shows the contours of $N_2(\widehat{\beta}, n(X'X)^{-1})$ (dashed red), the ``BIC prior''; note that the likelihood function (a function of $\beta$) is proportional to $N_p(\widehat{\beta}, (X'X)^{-1})$, so it has the same shape. When $r = -0.9$, the marginal likelihood of the true model is high with all the priors. If $r = 0.9$, the highest density regions of the likelihood of the true model are assigned relatively low probability density under $N_2(0_2, n(X'X)^{-1})$.

Table~\ref{tab:twopredictors} confirms this geometric intuition -- for sample sizes ranging from 5 to 15 and after 1000 simulations, the average posterior probability that the lower bound LB ($g$-prior with $g=n$) assigns to the true model is lower than with BIC  or the type II ML prior (ML). The Zellner-Siow (ZS) prior is less sensitive to the sign of $r$ than the lower bound, despite the fact that they are both centered at $0_p$ and have the same prior scale.

Our intuition can be supported mathematically. If $\sigma^2$ is known,
\begin{align} \label{eq:knowns2}
\log \left( \frac{\mathsf{BF}_{i0,\mathsf{BIC}}}{\mathsf{BF}_{i0,\mathsf{LB}}} \right) = \frac{p_i}{2} \log \left(\frac{n+1}{n} \right) + \frac{1}{2 \sigma^2} \left(1 - \frac{n}{n+1} \right) \mathsf{SSR}_{i} \ge 0,
\end{align}
which depends on the data only through $\mathsf{SSR}_i$. If the full model is true, $
\mathbb{E}[\mathsf{SSR}_i] = 2 \sigma^2  + n[\beta_1^2 + 2 \beta_1 \beta_2 r + \beta_2^2]$. Since in our example $\beta_1$ and $\beta_2$ are positive, $\mathbb{E}[\log (\mathsf{BF}_{i0,\mathsf{BIC}}/\mathsf{BF}_{i0,\mathsf{LB}})]$ increases as $r$ increases.
\end{example}

\begin{table}[h!]
\centering
\caption{Average posterior probability assigned to the true model (full model), $B = 1000$ simulations.}
\begin{tabular}{rrrrr|rrrr}
 &   \multicolumn{4}{c}{$r= -0.9$} & \multicolumn{4}{c}{$r= 0.9$} \\
  \hline
$n$ & BIC & ML & LB  & ZS & BIC & ML & LB  & ZS \\
  \hline
    5 & 0.954 & 0.797 & 0.665 & 0.503 & 0.978 & 0.911 & 0.361 & 0.310  \\
  10 &  0.997 & 0.983 & 0.976 & 0.952 & 0.999 & 0.995 & 0.605 & 0.964  \\
  15 &  1.000 & 0.999 & 0.999 & 0.997 & 1.000 & 1.000 & 0.874 & 0.998 \\
  20 & 1.000 & 1.000 & 1.000 & 1.000 &  1.000 & 1.000 & 0.979 & 1.000 \\
   \hline
   \hline
\end{tabular}
\label{tab:twopredictors}
\end{table}

The intuition we gathered from Example~\ref{ex:twopreds} that the type II ML procedure is between BIC and the lower bound (LB; i.e. a $g$-prior with $g=n$) is shown formally below.
\begin{proposition} \label{prop:compromise}
If $\mathcal{M}_j \supset \mathcal{M}_i$ (that is, if $\mathcal{M}_j$ contains all the predictors in $\mathcal{M}_i$), then
$$
\mathsf{BF}_{ji, \mathsf{BIC}} \ge \mathsf{BF}_{ji, \mathsf{ML}}\ge \mathsf{BF}_{ji, \mathsf{LB}},
$$
where  $\mathsf{BF}_{ji}= \mathsf{BF}_{j0}/\mathsf{BF}_{i0}$.  Let $\mathcal{M}_f$ be the full model (which includes all $p$ predictors) and $\mathcal{M}_0$ be the null model. If the prior on the model space is the same in all cases, the inequality above implies
\begin{align*}
\mathbb{P}_{\mathsf{BIC}}( \mathcal{M}_f \mid Y) \ge \mathbb{P}_{\mathsf{ML}}( \mathcal{M}_f \mid Y) \ge \mathbb{P}_{\mathsf{LB}}( \mathcal{M}_f \mid Y) \\
\mathbb{P}_{\mathsf{BIC}}( \mathcal{M}_0 \mid Y) \le \mathbb{P}_{\mathsf{ML}}( \mathcal{M}_0 \mid Y) \le \mathbb{P}_{\mathsf{LB}}( \mathcal{M}_0 \mid Y).
\end{align*}
\end{proposition}

If the true model is the full model, BIC assigns more probability to the truth than the type II ML prior and the lower bound; on the other hand, if the true model is the null model, the lower bound ($g$-prior) assigns more probability to the truth than the type II ML prior and BIC. However, there is yet another interesting asymmetry. When the true model is the null model, the differences between the lower bound and BIC tend to be small, whereas if the true model is the full model the differences can be rather large. We can provide some mathematical support to this claim. First, assume that $\sigma^2$ is known, so that the expression for $\log (\mathsf{BF}_{i0,\mathsf{BIC}}/\mathsf{BF}_{i0,\mathsf{LB}})$ is given in Equation~\ref{eq:knowns2}. If $\mu_\ast = X_0 \beta_{0\ast} + X_\ast \beta_\ast$ is the true value of the linear predictor, we can write $\mathbb{E}[\mathsf{SSR}_i] = p_i \sigma^2 + \lVert \mathsf{P}_{{X}_i} \mu^\ast \rVert^2$. For fixed $X_i$, $\mathbb{E}[\log (\mathsf{BF}_{i0,\mathsf{BIC}}/\mathsf{BF}_{i0,\mathsf{LB}})]$ is minimized when $\beta_{\ast} = 0_{p_\ast}$ (i.e. when the null model is true). Also note that $\mathbb{E}[\mathsf{SSR}_i]$ is increasing in $p_i$, which implies that the expected (log) differences between the lower bound and BIC grow as the number of predictors grows. For unknown $\sigma^2$, $\log (\mathsf{BF}_{i0,\mathsf{BIC}}/\mathsf{BF}_{i0,\mathsf{LB}})$
is increasing in $R^2_i$, which is consistent with our argument.

At the beginning of this section, we mentioned that a type II ML prior that has been previously studied is the $g$-prior $N_p(\beta \mid 0_p, \widehat{g} \,  \sigma^2 (X'X)^{-1})$, where $\widehat{g}$ is set locally by maximizing the marginal likelihood subject to $g \ge 0$ \citep{george2000calibration, hansen2003minimum, liang2008mixtures}. This prior has undesirable features that are a byproduct of not maximizing the marginal likelihood subject to a lower bound on $g$ that is bounded away from 0. One of them is that the resulting null-based Bayes factors are always greater or equal to 1 (which leads to inconsistency if the null model is true), and another one is that the Bayes factor between any two models can be equal to 1 with positive probability in cases where $n > p + p_0$ (especially when $n \approx p + p_0$), which cannot occur (with positive probability) with proper priors or our restricted type II prior.

We close this subsection by studying whether the type II ML prior satisfies the desiderata in \cite{bayarri2012criteria} for objective priors in model selection.
\begin{enumerate}
\item \textbf{Basic criterion:} The basic criterion is satisfied if the prior is proper, which the type II ML prior satisfies directly because of the restriction.
\item \textbf{Model selection consistency:} Let the true model be $\mathcal{M}_\ast: N_n( Y \mid X_0 \beta_0 + X_\ast \beta_\ast, \sigma^2 I_n)$. Then, model selection consistency is satisfied if $\mathbb{P}(\mathcal{M}_\ast \mid Y)$ converges to 1 in probability. The type II ML prior is model-selection consistent under the following regularity condition, which is commonly made in the literature \citep{fernandez2001benchmark,liang2008mixtures,guo2009bayes,maruyama2011fully,bayarri2012criteria,som2014block}. For any model $\mathcal{M}_j$ that doesn't nest the true model, assume that
$$
\lim_{n \rightarrow \infty} \frac{\beta_\ast' X_\ast' (I_n-\mathsf{P}_{X_j})X_\ast \beta_\ast}{n} = b_j \in (0, \infty).
$$
The assumption can be interpreted as that the models have design matrices that can be differentiated in the limit \citep{bayarri2012criteria}.
\item \textbf{Information consistency:} Suppose that, for a fixed $n$, $\lVert \widehat{\beta}_i \rVert \rightarrow \infty$, which implies $R^2_i \rightarrow 1$. This is a situation where there is overwhelming evidence in favor of $\mathcal{M}_i$ \citep{liang2008mixtures}. Information consistency holds if $\mathsf{BF}_{i0} \rightarrow \infty$, which is satisfied by the type II ML prior.
\item \textbf{Intrinsic consistency:} A prior satisfies intrinsic consistency if, as $n$ grows, it converges to a proper prior which does not depend on model-specific parameters or $n$. In general, this criterion isn't satisfied for the type II ML prior. To see this, assume that $(X_i'X_i)/n \rightarrow \Xi_i$ for a positive definite matrix $\Xi_i$, which holds if there is a fixed design or the covariates are drawn independently from a distribution with finite second moments \citep{bayarri2012criteria}. Then, the prior covariance $\widehat{W}_\ast$ for the true model has the limiting behavior (in probability)
$$
\widehat{W}_\ast \rightarrow_P \begin{cases}  \Xi_\ast^{-1} & \text{if } \, \beta_\ast'\Xi_\ast \beta_\ast \le \sigma^2_\ast \\ \left( \frac{1}{\sigma^2_\ast} - \frac{1}{\beta_\ast' \Xi_\ast \beta_\ast} \right) \beta_\ast \beta_\ast' + \Xi_\ast^{-1} & \text{if } \, \beta_\ast'\Xi_\ast \beta_\ast > \sigma^2_\ast \end{cases},
$$
which depends on $\beta_\ast$ and $\sigma^2_\ast$.
\item \textbf{Null and dimensional predictive matching:} In both cases, the notion of minimal training sample size is central to the definition. For any model $\mathcal{M}_i$, the minimal training sample size is the smallest sample size $n^\ast_i$ such that the marginal likelihood of the model is finite. Null predictive matching is achieved if, for any model $\mathcal{M}_i$, we have $\mathsf{BF}_{i0} = 1$ when the sample size is equal to the minimal training sample size $n^\ast_i$. Dimensional predictive matching is achieved if, for any pair models of the same dimension $\mathcal{M}_i$ and $\mathcal{M}_j$, we have $\mathsf{BF}_{ij} = 1$ whenever $n^\ast_i = n^\ast_j$. The type II ML prior isn't null or dimensional predictive matching. For $p > 1$, the minimal training sample size for the type II ML prior is $n = p +p_0+1$. [If $p=1$, the marginal likelihood doesn't depend on the choice of $W$.] When $n = p+p_0$, the marginal likelihood is finite for any given $W$, but one can choose $W \succeq n (X'X)^{-1}$ so that the marginal goes to $\infty$ (this is shown in the supplementary material). Null predictive matching isn't satisfied: in fact, $\mathsf{BF}_{i0}$ goes to $\infty$ as $R^2_i \rightarrow 1$ when $n = p + p_0 +1$. Similarly, it is easy to see that dimensional predictive matching isn't satisfied, either; different models will have different $R^2_i$, yielding Bayes factors that are different than 1.
\item \textbf{Invariance:} The type II ML prior is invariant with respect to linear transformations of the design matrix (e.g. changes of measurement units). More explicitly, let $A$ be an invertible $p \times p$ matrix and $\tilde{X} = X A$. Let $\beta$ and $\tilde{\beta}$ be the regression coefficients of the linear model if the design matrices are $X$ and $\tilde{X}$, respectively. If the type II ML prior is put on $\beta$ and $\tilde{\beta}$, then $\beta$ and $A \tilde{\beta}$ are equal in distribution.
\end{enumerate}

 Table~\ref{tab:properties} compares the properties of the type II ML prior with those of BIC, the lower bound LB ($g$-prior with $g=n$), the Zellner-Siow prior (ZS), and the type II ML $g$-prior where $g$ is set locally by maximizing the marginal likelihood subject to $g \ge 0$, which we denote $\widehat{g}$. Our type II ML prior is model-selection consistent, whereas the $\widehat{g}$-prior isn't under the null model; however, the $\widehat{g}$-prior is predictive matching, while our type II ML prior isn't. According to the definition above, it doesn't make sense to assert that BIC is invariant to linear transformation (since it isn't a prior), but it depends on the data only through $R^2$, which is invariant with respect to invertible linear transformations.

 \begin{table}
\centering
\caption{Comparison of model selection desiderata for different approaches.}
    \begin{tabular}{l|lllll}
    ~                       & ML & BIC            &  LB & ZS  &  $\widehat{g}$                 \\
\hline
    Proper                  & yes                         & -            & yes   & yes & yes                  \\
    Model selection consistency   & yes                         & yes            & yes   & yes & no                   \\
    Information consistency & yes                         & yes            & no    & yes & yes                  \\
    Intrinsic consistency   & no         & -             & yes   & yes & no  \\
    Predictive matching     & no                          & no             & yes   & yes & yes                   \\
   	Invariance  			& yes                         & -            & yes   & yes & yes                  \\
    Closed form Bayes factors  & yes    & yes               &yes & no  &yes \\
    \hline
    \end{tabular}
    \label{tab:properties}
\end{table}

It is not a surprise that data-dependent priors lack some of the desirable properties of real priors. One sacrifices some Bayesian
 features when leaving the pure Bayesian domain.

\subsection{Estimation and prediction} \label{sec:estimation}

\subsubsection{The type II ML posterior mean}

For simplicity, we omit model subscripts and assume that the model is $Y \sim N_n(X_0 \beta_0 + X \beta, \sigma^2 I_n)$, $X_0'X = 0_{p_0 \times p}$. If we put the right-Haar prior $\pi(\beta_0, \sigma^2) \propto 1/\sigma^2$ on the common parameters and the type II ML prior on $\beta \mid \sigma^2$, the posterior mean of $\beta$ is
$$
\tilde{\beta} = \mathbb{E}(\beta \mid Y) = \begin{cases}    \frac{n}{n+1} \widehat{\beta} & \text{if } R^2 \le \frac{n+1}{2n-p_0} \\ \left(1-\frac{1-R^2}{(n-p_0-1)R^2}\right) \widehat{\beta} & \text{if } R^2 > \frac{n+1}{2n-p_0}. \end{cases}
$$
The expression can be derived by applying the Sherman-Morrison formula twice to $\mathbb{E}(\beta \mid Y) = [\widehat{W}^{-1}+X'X]^{-1} X' Y$. The properties of an analogous estimator in the normal means problem (for known $\sigma^2$) are studied in \cite{dasgupta1988frequentist}, where it is shown that it is minimax with respect to squared error loss. Proposition~\ref{prop:minimax} shows that $\mathbb{E}(\beta \mid Y)$ is also minimax with respect to a (scaled) predictive loss because it belongs to the class of minimax estimators characterized in \cite{strawderman1973proper}.

\begin{proposition} \label{prop:minimax}
Let $p \ge 3$ and $n > p + p_0$. The estimator $\tilde{\beta} = \mathbb{E}(\beta \mid Y)$ is minimax with respect to the (scaled) squared predictive loss
$$
L(\beta, \delta) = (\beta-\delta)' (X'X) (\beta - \delta)/\sigma^2 \,.
$$
\end{proposition}

 The mean squared error of the posterior mean of the lower bound prior (Zellner's $g$-prior, where $g = n$) is increasing in $\lVert \beta \rVert$. On the other hand, the mean squared error of $\widehat{\beta}$ is constant in $\lVert \beta \rVert$. The estimator $\tilde{\beta}$ is equal to the posterior mean of the lower bound when $R^2$ is small, and close to $\widehat{\beta}$ when $R^2$ is large. Therefore, $\tilde{\beta}$ avoids ``selecting'' the lower bound in cases where it has high mean squared error (that is, whenever $\lVert \beta \rVert$ and $R^2$ are large).

 \subsubsection{A simulation study with correlated predictors} \label{sec:simstudies}

To gain further insight into the differences between the type II ML prior, the lower bound (LB) prior ($g$-prior with $g = n$), the Zellner Siow (ZS) prior and BIC,
we simulate data from $Y = 1_n \, \alpha + X \beta + \epsilon$, $\epsilon \sim N_n(0_n, \sigma^2 I_n)$, where $n = 50$, $\alpha = 2$, $\sigma^2 =1$, and $\beta$ is $8$-dimensional with $k$ nonzero elements, for $k \in \{0,1,2, \, ... \, , 8\}$. We consider 2 different types of correlation between the predictors: the orthogonal case $X'X = I_p$ and an AR(1) structure
$$
\frac{1}{n-1}(X'X) = \left( \begin{array}{ccccc}
1 & \rho & \rho^2 & ... & \rho^p \\
\rho & 1 & \rho & ... & \rho^{p-1} \\
\rho^2 & \rho & 1 &  ... & \rho^{p-2} \\
\vdots & \vdots & \vdots & \ddots & \vdots \\
\rho^p & \rho^{p-1} & \rho^{p-2} & ... & 1 \end{array} \right)
$$
for $\rho = 0.9$. For all $k$, we generate $\beta_k \sim N_k(0_k, g I_k)$. The location of the $k$ zeros in the $\beta$ vector is drawn at random (according to the uniform distribution). We use $g \in \{5,25\}$ as in \cite{cui2008empirical} and \cite{liang2008mixtures}, representing weak and strong signal-to-noise ratios, and evaluate performance with respect to the predictive squared loss function
$
L(\beta, \delta) = \lVert X \beta - X \delta \rVert^2,
$
where $\delta$ is an estimator of $\beta$. [This is also the loss function that was used in the simulation studies in \cite{cui2008empirical} and \cite{liang2008mixtures}]. The estimators that are considered for the various priors are the posterior means (and $\widehat{\beta}$ in the case of BIC) of the highest probability model (HPM) and the median probability model (MPM), and the estimate arising from Bayesian model averaging (BMA).
We ran $1000$ simulations for all scenarios and the results are displayed in Figures 1 and 2 in the supplementary material.

In the orthogonal case, BIC, the type II ML prior, LB ($g$-prior with $g=n$) and ZS behave similarly when $g = 5$. When $g = 25$, we can observe more differences: LB is progressively worse than the rest as the number of true predictors increases, ZS is slightly better than BIC and the type II ML prior when not all predictors are active, and the difference between ZS and BIC and the type II ML prior narrows as the number of true predictors increases.

 The results with the AR(1) correlation structure show bigger discrepancies. As the number of true predictors increases, the loss of the LB is substantially higher than the loss with any other prior, especially when $g=25$. When $g = 5$, both LB and ZS are outperformed by BIC and the type II ML prior. When $g = 25$, ZS has similar losses as BIC and the type II ML prior when the number of true predictors is between 0 and 6, but is outperformed when the true number of predictors is 7 or 8 (in which case, the true model is the full model).

 In the cases where the LB is clearly outperformed, its posterior distribution over the model space is closer to the uniform distribution than the other posteriors, as evidenced in the first panel in Figure 3 in the supplementary material, which shows the average entropy of the posterior distributions over the model space. Additionally, ZS induces a noticeably less entropic (more concentrated) posterior distribution over the model space, especially when few predictors are active.  ZS and the LB select HPMs and MPMs with fewer predictors than BIC and the type II ML prior (see second and third panel in Figure 3 in the supplementary material, which show the percentage of times the MPM equals the true model and the average size of the MPM, respectively). When the true model is the full model, an interesting phenomenon occurs: ZS is the prior where the MPM is equal to the true model less often, but the average predictive loss of the prior stays competitive with BIC and ML. Upon further inspection in our simulations, this is due to the fact that when some of the true coefficients are non-zero but rather small, ZS does not include their predictors in its MPM, but that does not worsen the predictive loss by much. The HPM and MPM with BIC and the type II ML prior tend to be the same model, and they coincide with the models selected with the LB in the cases where the signal is low, as expected. On the other hand, when the signal is high, the LB assigns more probability to wrong models than the other approaches, and sometimes the HPM and MPM end up being an egregiously bad model, resulting in a substantially higher average loss. Note that ZS, which also has $n(X'X)^{-1}$ as its prior scale but has thicker tails, does not seem to be nearly as affected by this issue as the LB, especially when the signal is high enough (i.e. $g=25$).

\section{High-dimensional ANOVA} \label{sec:ANOVA}

In this section, we revisit the one-way ANOVA problem that was introduced in \cite{stone1979comments} and later studied in \cite{berger2003approximations}. In this example, the number of predictors $p$ grows to infinity. Suppose we have observations
$$
y_{ij} = \mu_i + \epsilon_{ij}, \, \epsilon_{ij} \stackrel{\mathrm{iid}}{\sim} N_1(0,1)
$$
where $i \in \{1,2, \, ... \, , p\}$ (groups) and $j \in \{1,2, \, ... \, , r\}$ (replicates). We assume that $r$ is fixed and $p$ grows to infinity. We only consider the null model $\mathcal{M}_1 : \mu = 0_p$ and the full model $\mathcal{M}_2 : \mu \neq 0_p$. If the true model is $\mathcal{M}_2$, we assume that $\mathrm{lim}_{p \rightarrow \infty} \lVert{\mu}\rVert^2/p = \tau^2 > 0.$

Let $\ell$ be the log-likelihood function of the full model and $\widehat{\mu}$ the maximum likelihood estimate of $\mu$. If BIC is defined as $-2 \ell(\widehat{\mu}) +p \log n$,
it is inconsistent under $\mathcal{M}_2$ \citep{stone1979comments}. \cite{berger2003approximations} show that, if the prior on $\mu$ under $\mathcal{M}_2$ is $\mu \mid g \sim N_p(0, g I_p)$ with a mixing density over $g$ (which doesn't depend on $n$) with support $(0, +\infty)$, consistency holds. Alternatively, if $g$ has restricted support $(0, T)$ for $T < \infty$, there is a region of inconsistency under $\mathcal{M}_2$.

In this problem, the prior scale has to be chosen carefully. A naive parallel of the type II ML prior in Section~\ref{sec:linmod} would have as lower bound for the prior covariance $n (X'X)^{-1} = ({n}/{r}) I_p = p I_p.$
However, it is straightforward to show that any normal prior whose scale goes to infinity as $p \rightarrow \infty$ is inconsistent under $\mathcal{M}_2$. Since the effective sample size of $\mu$ in this problem is $r$ instead of $n$ (see \cite{berger2014effective}), we take $g = r$ and study the properties of a prior whose covariance is $r (X'X)^{-1} = I_p.$ In the same vein, BIC can be defined appropriately by taking $\log r$ as the penalty instead of $\log n$. The asymptotic behavior of both approaches can be summarized as follows:
\begin{itemize}
	\item Normal prior with $I_p$ as prior covariance: Under $\mathcal{M}_1$, consistency for all $r \ge 1$. Under $\mathcal{M}_2$, inconsistency if $\tau^2 \le (1+r) \log (1+r)/r^2 -1/r$ and consistency otherwise. For example, if $\tau^2 = 0.25$, consistency holds under $\mathcal{M}_2$ for $r \ge 5$, and consistency holds for all $r$ if $\tau^2 > 2 \log 2 -1$.
	\item BIC with $\log r$ as penalty: Under $\mathcal{M}_1$, inconsistency if $r \in \{1,2\}$ and consistency otherwise. Under $\mathcal{M}_2$, inconsistency if $\tau^2 \le (\log r -1)/r$ and consistency otherwise. The condition is most stringent at $r = e^2$, so consistency holds for all $r$ if $\tau^2 > 1/e^2$.
\end{itemize}
Under $\mathcal{M}_2$, the region of inconsistency of BIC is contained in the region of inconsistency of the normal prior; however, BIC can be inconsistent under $\mathcal{M}_1$. 

The type II ML prior
\begin{align*}
\mu &\sim N_p( 0_p, \widehat{W}) \\
\widehat{W} &= \mathrm{arg \, max}_{W \succeq I_p} m(Y) = I_p + \max\{ 0, 1-(r+1)/(r\lVert \widehat{\mu} \rVert^2) \} \widehat{\mu} \widehat{\mu}'
\end{align*}
yields the Bayes factor
$$
\mathsf{BF}_{21} = \begin{cases}
(r+1)^{-p/2} \exp \left\{ \frac{r^2 \lVert \widehat{\mu} \rVert^2 }{2(r+1)} \right\}  & \text{if } \lVert \widehat{\mu} \rVert^2 \le  1 + 1/r  \\
\left( \frac{r \lVert \widehat{\mu} \rVert^2}{r+1} \right)^{-1/2} (r+1)^{-p/2} \exp \left\{ \frac{(r \lVert \widehat{\mu} \rVert^2 - 1)}{2} \right\}  & \text{if } \lVert \widehat{\mu} \rVert^2 > 1 + 1/r
\end{cases} \,.
$$
Under $\mathcal{M}_1$, the type II ML Bayes factor is inconsistent for $r = 1$ and consistent for all $r > 1$. Under $\mathcal{M}_2$, it is inconsistent for $\tau^2 \le [\log(r+1)-1]/r$ and consistent otherwise.

The type II ML prior acts as a compromise between the normal prior and BIC but, unfortunately, it still has regions of inconsistency which mixtures of normal priors avoid. However, the type II ML Bayes factor is available in closed form, whereas the Bayes factors that stem from using mixtures of normals generally are not. 

PBIC and PBIC*, which are prior-based versions of BIC that are defined and studied in \cite{bayarri2019prior}, are consistent under $\mathcal{M}_1$ for all $r \ge 1$, but inconsistent under $\mathcal{M}_2$ for $\tau^2 < [\log 2 + \log(r+1)-1]/r$. That is, under $\mathcal{M}_2$, the region of inconsistency of PBIC and PBIC* contains the region of inconsistency of the restricted type II ML prior. On the other hand, our type II ML prior is inconsistent under $\mathcal{M}_1$ for $r = 1$, while PBIC and PBIC* are not. Therefore, in this example, the type II ML prior discussed here is more favorable to $\mathcal{M}_2$ than PBIC and PBIC*.


\section{Incorporating prior information} \label{sec:shibata}

The constraints we have placed on the type II ML prior have been basic constraints, preventing the prior from becoming too concentrated.
It is also possible to use constraints that incorporate available prior information, which can lead to improved inferences. We illustrate this possibility by revisiting
the example in \cite{shibata1983asymptotic}, which was also studied in \cite{barbieri2004optimal}.

The goal in the Shibata example is to estimate the function $f(x) = - \log(1-x)$, $-1 \le x \le 1$ from independent observations $y_i = f(x_i) + \varepsilon_i$, where the $\varepsilon_i$ are independent  $\varepsilon_i \sim N_1(0, \sigma^2 I_n)$ and $\sigma^2$ is known. The function $f$ can be expressed in an orthogonal series expansion as $
f(x) = \sum_{i = 1}^\infty \beta_i \phi_i(x),
$
where $\phi_i(x)$ are the Chebyshev polynomials of the first kind. We approximate $f$ with a finite series expansion, modeling $y_i = \sum_{i = 1}^j \beta_i \phi_i(x) + \varepsilon_i$. We consider different truncation points $j$, ranging from $1$ to $k$, so our model space consists of a sequence of nested models
$$
\mathcal{M}_j : Y \mid \alpha, \beta_j, \sigma^2 \sim N_n(1_n \alpha + X_j \beta_j, \sigma^2 I_n)
$$
for $j \in \{1,2,3, \, ... \,  , k\}$, where the design matrices $X_j$ have dimension $n \times j$ and the columns are given by the Chebyshev polynomials of the first kind evaluated at the knots $x_i = \cos ( \pi (n-i+1/2)/n)$, for $i \in \{1,2, \, ... , n\}$. The true coefficients in an infinite orthogonal expansion are $\alpha = \log 2$ and $\beta_j = 2/j$. The design matrices are orthogonal with $X_j ' X_j = (n/2) I_{j}$ and $1_n' X_j = 0_j'$. [See \cite{barbieri2004optimal} for a more detailed explanation.]

We consider $n = 30, k=29, \sigma^2 = 1$, $n = 100, k=79, \sigma^2 = 1$, and $n = 2000, k = 79, \sigma^2 = 3$ and put a uniform prior (i.e., 1/29 or 1/79) on the size of the nested models. We utilize two local type II ML priors based on $\beta \sim N(0, \sigma^2 A)$:
\begin{itemize}
\item The unit-information constraint $A \succeq  n (X'X)^{-1}$.
\item In polynomial regression, the true coefficients often decrease at polynomial rate. With that in mind, we define a type II ML prior whose covariance matrix is diagonal, with diagonal elements decreasing according to some power law. That is, $A = \mathrm{diag}(d_1, d_2, \, ... \, , d_k)$ with $d_i = c i^{-a}$ for $i \in \{1, 2, \, ... \, , k\}$. The parameters $c, a \ge 0$ are found by maximizing the marginal likelihood.
\end{itemize}
We will compare these three methods on Shibata's example, utilizing squared predictive loss $L(f, \widehat{f}) = \int_{-1}^{1} (f(x) - \widehat{f}(x))^2 \, \mathrm{d} x$, as in \cite{barbieri2004optimal}.
We also consider AIC ($-2\ell(\widehat{\beta}_j) + j$) and BIC ($-2\ell(\widehat{\beta}_j) +j \log n$), treating $\exp(-\mathrm{AIC}/2)$ and $\exp(-\mathrm{BIC}/2)$ as approximate marginal likelihoods. We compare the predictive loss of Bayesian model averaging (BMA), the median probability model (MPM; \cite{barbieri2004optimal}), and the highest probability model (HPM). Note that the AIC and BIC columns for the HPM
correspond to use of the actual AIC and BIC criteria, since maximizing the posterior probability is
equivalent to minimizing the criterion. The MPM and BMA columns utilize AIC and BIC by converting them to approximate
marginal likelihoods and utilizing the relevant Bayesian theory.

The results are summarized in Table~\ref{tab:shibata}. BIC and $\widehat{W}$ behave similarly in all cases, as we have seen in previous sections. The informative type II ML priors outperform the others. AIC is somewhat better than BIC, and their Bayesian implementations (MPM and BMA)
outperform use of the raw criteria (HPM).

 All across the board, BMA outperforms the rest (as expected), followed by the MPM and the HPM; the MPM
 is the best single predictive model in nested model scenarios, as shown in \cite{barbieri2004optimal}.

\begin{table}[h!]
\centering
\caption{Predictive loss, based on $N = 1000$ simulations. Average model sizes in square brackets.}
\label{tab:shibata}
\begin{tabular}{l|lllll}
    \textbf{HPM}  & $c i^{-a}$ & $\widehat{W}$   & AIC  & BIC  \\
\hline
$n = 30, k = 29, \sigma^2 = 1$  & \textbf{0.904} [10]  & 1.141 [4] & 1.076 [7] & 1.131 [4] \\
$n = 100, k = 79, \sigma^2 = 1$  & \textbf{0.471} [23] & 0.693 [7]  & 0.582 [13] & 0.692 [7]  \\
$n = 2000, k =79, \sigma^2 = 3$ & \textbf{0.136} [57] & 0.295 [13]  & 0.188 [36] & 0.295 [13]  \\
\textbf{MPM}    &   &   &     &   &   \\
\hline
$n = 30, k = 29, \sigma^2 = 1$  &  \textbf{0.839} [16] & 1.093 [4] & 1.027 [7] & 1.089 [4] \\
$n = 100, k = 79, \sigma^2 = 1$   & \textbf{0.441} [44] & 0.680 [7]  & 0.566 [13] & 0.679 [7]  \\
$n = 2000, k =79, \sigma^2 = 3$   & \textbf{0.134} [59] & 0.289 [13]  & 0.185 [37] & 0.289 [13]  \\
\textbf{BMA}   &  &  &    &  & \\
\hline
$n = 30, k = 29, \sigma^2 = 1$   & \textbf{0.837} & 0.990  & 0.921 & 0.983  \\
$n = 100, k = 79, \sigma^2 = 1$  & \textbf{0.437} & 0.623  & 0.521 & 0.621 \\
$n = 2000, k = 79, \sigma^2 = 3$   & \textbf{0.133}  & 0.275   & 0.170 & 0.275   \\
\end{tabular}
\end{table}

\section{Conclusions} \label{sec:conclusions}

Conceptually, the type II ML priors we studied offer an attractive compromise between conventional priors, which might seem overly concentrated at the null model, and BIC. The importance of constraining
the maximization so that the prior does not overly concentrate was highlighted, and the need to carefully choose the constraint
in high-dimensional situations was discussed.

The surprise of the analysis was that the type II ML prior gives remarkably similar answers to BIC. Indeed,
the paper could be viewed as primarily providing a new justification of BIC in normal linear
models, suggesting that BIC need not just be viewed as an approximation but as something that corresponds quite
closely to an interpretable type II ML procedure (and not just with priors that sit on top of the model likelihoods).

In Example~\ref{ex:twopreds} and the simulation study in Section~\ref{sec:simstudies}, we observe that the $g$-prior with $g=n$, which is the lower bound of our restricted type II ML prior, can severely underperform when the predictors are correlated (especially when most predictors are active
). In our numerical comparisons, the Zellner-Siow prior, BIC, and the type II ML procedure yield similar results. Zellner-Siow seems to perform slightly better in most cases, but its performance suffers when most predictors are active. 
From a theoretical perspective, Zellner-Siow satisfies intrinsic consistency and predictive matching, which are not satisfied by the type II ML prior. However, the type II ML prior yields closed form Bayes factor, whereas the Zellner-Siow prior does not (see Table~\ref{tab:properties}).

Finally, we revisited the nonparametric regression example in \cite{shibata1983asymptotic}, showing how
prior information could be incorporated into the constraints defining type II ML priors, leading
to considerably improved performance (when the prior information is correct). This is perhaps the most
promising practical venue for type II ML priors: embed available structural information about the prior
into the class of priors, and then use type II ML.

\bibliographystyle{chicago}

\bibliography{restrictedML}

\end{document}